\documentclass{article}
\usepackage{amsbsy,amssymb,amscd,amsfonts,latexsym,amstext,delarray,
amsmath}  \headsep=15pt
\def\hpi{\hat{\pi}}

\def\bX{ \bar{X}}

\def\Spec{{ \mbox{Spec} }}

\def\ra{{ \rightarrow }}

\def\d{{ \delta }}

\def\F{ {\bf F} }
\def\hZ{ \hat{\Z}}

\def\G{{ \Gamma }}
\def\Gal{{ \mbox{Gal} }}

\def\bQ{\bar{\Q}}
\def\bF{{ \bar{F} }}

\def\Z{{ \mathbb{Z}}}

\usepackage{amsmath}
\usepackage{amsfonts}
\usepackage{amscd}
\usepackage{amssymb}
\def\Q{\mathbb{Q}}

\def\invlim{\varprojlim}

\def\bF{ \bar{F}}

\title{ Diophantine geometry as Galois theory in the mathematics of Serge Lang}
\author{ Minhyong Kim}

\begin{document}
\maketitle
Lang's conception of Diophantine geometry is rather compactly represented
by the following celebrated conjecture
\cite{DP}:
\begin{quotation}
Let $V$ be a subvariety of a semi-abelian variety $A$, $G \subset A$ a finitely generated
subgroup, and $Div(G)$ the subgroup of $A$ consisting of the division points of $G$. Then $V\cap Div(G)$ is
contained in a finite union of subvarieties of $V$ of the form $B_i+x_i$, where each $B_i$ is a semi-abelian
subvariety of $A$  and $x_i\in A$.
\end{quotation}
There is  a wealth of literature at this point surveying the various ideas and techniques employed in its
resolution, making it unnecessary to review them here in any detail \cite{faltings3, mcquillan}.
However, it is still
worth taking note of the valuable {\em generality } of the formulation, evidently arising  from a profound instinct
for the plausible structures of mathematics. To this end, we remark merely that it was exactly this
generality that made possible the astounding interaction with geometric model theory in the 90's \cite{bouscaren}. That is to
say, analogies to model-theoretic conjectures and structure theorems would have been far harder to detect if
attention were restricted, for example, to situations where the intersection is expected to be finite.
Nevertheless, in view of the sparse subset of the complex net of ideas surrounding this conjecture that we
wish to highlight in the present article, our intention is to focus exactly on the case where $A$ is compact
and $V$ does {\em not} contain any translate of a connected proper subgroup. The motivating example, of course, is a
compact hyperbolic curve embedded in its Jacobian. Compare then the two simple cases of the conjecture that
are amalgamated into the general formulation:
\begin{quotation}
(1) $V\cap A[\infty]$, the intersection between
 $V$ and the torsion points of $A$,
is finite.

(2) $V\cap G$ is finite.
\end{quotation}
 Lang expected conjecture (1) to be resolved
using  Galois theory alone. This insight was based upon work of Ihara, Serre, and Tate (\cite{FDG},
VIII.6),
dealing with the analogous problem for a torus, and comes down to the
conjecture, still unresolved,
 that the image of the Galois
representation in $Aut (A[\infty])\simeq GL_{2g}(\hat{\Z})$ contains an open subgroup of the homotheties
$\hat{{\Z}}^*$. Even while assertion (1) is already
a theorem of  Raynaud \cite{raynaud},
 significant progress along the lines originally
envisioned by Lang was
 made in \cite{bogomolov} by
 replacing $A[\infty]$ with $A[p^{\infty}]$, the points of $p$-power torsion,
and  making crucial use of
$p$-adic Hodge theory.

It is perhaps useful to reflect briefly on the overall
context of Galois-theoretic methods in Diophantine geometry, of course without attempting to do justice to the
full range of interactions and implications. Initially, that Galois theory is  relevant to the study of
Diophantine problems should surprise no one. After all, if we are interested in
 $X(F)$, the set of rational points of a variety
$X$ over a number field $F$, what is more natural than to observe that $X(F)$ is merely the fixed point set  of
$\G:= Gal(\bar{F}/F)$ acting on $X(\bar{F} )$? Since the latter is an object of classical geometry,
such an expression might be expected to give us a tight grip on the subset $X(F)$ of interest. This view is of
course very naive and the action on $\G$ on $X(\bF)$ is notoriously difficult to use in any direct fashion.
 The action on {\em torsion points} of
commutative group varieties on the other hand, while still difficult, is considerably more tractable than the
general situation, partly because the added structure of a finite abelian group
behaves well under specialization. There is also serious input from topology in that the action on torsion points
can also be viewed as an action on
 homology groups
$H_1(\bX, \Z/m)$ with finite coefficients, or their inverse limits
$H_1(\bX, \hZ):=\invlim H_1(\bX, \Z/m)$. It
is by now known even to the general public that an intricate study of this action underlies the theorem of Wiles \cite{wiles}, and
roughly one-half of the difficulties in the theorem of Faltings
\cite{faltings1}. Within
such a context, it is understandable that
Galois theory might be conceived of as exceedingly powerful a tool for (1).

On the other hand, for conjecture (2), where the points to be studied are not torsion, it is not at all clear that Galois theory
can be as useful. In fact, my impression is that Lang expected {\em analytic geometry} of
some sort to be the main input to conjectures of type (2). This is indicated, for example, by the absence of any
reference to arithmetic  in the statement. We could say that implicit in the conjecture is
a powerful  idea that we may refer to as the {\em analytic strategy}:
\begin{quote}
(a) replace the difficult Diophantine set
$V(F)$ by the geometric intersection $V\cap G$; \\
(b) try to prove this intersection finite by analytic means.
\end{quote}
In this form, the strategy appears to have been extraordinarily efficient over function fields, as in the work
of Buium \cite{buium}. Even Hrushovski's \cite{hrushovski}
proof can interpreted in a similar light where passage to the completion of
suitable theories is analogous to the move from algebra to analysis
(since the theory of  fields is not good
enough). These examples should already
suffice to convince us  that it is best left open as to what kind of analytic means are most appropriate
in  a given situation. The proof  over number fields by Faltings \cite{faltings2},  as well as the curve
case by Vojta \cite{vojta}, employ rather heavy Archimedean analytic geometry. Naturally, the work of Vojta and Faltings
draws us away from the realm of traditional
Galois theory. However, in
Chabauty's theorem \cite{chabauty},
 where $V$ is a curve and the rank of $G$ is strictly less than the dimension of $A$, it is elementary
non-Archimedean analysis, more specifically
$p$-adic abelian integrals,  that completes the proof. Lang makes clear in several different places
(\cite{IP}, \cite{FDG}, notes to chapter 8, \cite{NT}, I.6) that Chabauty's theorem
was a definite factor
 in the formulation of his conjecture. This then invites a return
to our main theme,  as we
remind ourselves
that
 non-Archimedean analysis has come to be viewed profitably over the last several decades as
a projection of {\em analysis on Galois groups}, a perspective of which Lang was well aware
(\cite{CF}, chapter 4). As such, it  has something quite substantial to say about
non-torsion points, at least on elliptic curves (\cite{kato}, for example).
Hodge theory is again a key ingredient, this time as the medium in which to
realize such a projection \cite{perrin-riou}.

The  first step in the Galois-theoretic
description of non-torsion points,
at once elementary and fundamental, goes through the Kummer exact sequence
$$0\ra A[m](F) \ra A(F) \ra A(F) \stackrel{\d}{\ra} H^1(\G, A[m]) \ra H^1(\G,A)\ra .$$
 In this case, an easy study of specialization allows us to locate the image of $\d$ inside a subgroup
$H^1(\G_S, A[m])$ of cohomology classes with restricted ramification, which then form a finite group. We deduce thereby
the finiteness of $A(F)/mA(F)$, the weak Mordell-Weil theorem. Apparently, a streamlined presentation of this proof,
systematically emphasizing the role of Galois cohomology,  first appears in Lang's paper with Tate
\cite{LT}. There,
they also emphasize the interpretation of Galois cohomology groups as classifying spaces for {\em torsors}, in
this case, for $A$ and $A[m]$. (We recall that a torsor for a group $U$ in some category is an object
corresponding to a set with simply transitive $U$-action, where the extra structure of the category, for
example, Galois actions, prevent them from being trivial. (\cite{milne}, III.4)) In the inverse limit indexed by $m$ running over
powers of a fixed prime $p$, and incorporating the information from various towers of field extensions, we
encounter the analytic arguments that yield deep results for elliptic curves.

It is perhaps not entirely obvious that an equally elementary
analogue of the Kummer map $\d$ exists for hyperbolic curves
that is furthermore expected to provide yet another instance of the
analytic strategy.
In the course of preparing this article, I looked into Lang's
{\em magnum opus} \cite{FDG} for the first time in many years and was
a bit surprised to find a section entitled `non-abelian Kummer theory.' What is non-abelian there is the Galois
group that
 needs to be considered if one does not assume a priori
that the torsion points of the group variety are  rational over the ground field. The field of $m$-divison
points of the rational points will then have a Galois group $H$ of the form
$$0\ra A[m] \ra H \ra M \ra 0$$
where $M\subset GL_{2g}(\Z/m)$. Thus, `non-abelian' in this context is used in the same sense as in the reference to
non-abelian Iwasawa theory. But what is necessary for hyperbolic curves is yet another layer of
non-commutativity, this time in the coefficients of the action. Given a variety $X$ with a
rational point $b$, we can certainly consider the \'etale fundamental group $\hpi_1(\bX,b)$ classifying  finite
\'etale covers of $\bX$. But the essential construction is that given any other point $x\in X(F)$, we end up with the set of
\'etale paths $$\hpi_1(X;b,x)$$ from $b$ to $x$ which is naturally a torsor for $\hpi_1(\bX,b)$. All these live
inside the category of pro-finite sets with Galois action. There is then a non-abelian continuous
cohomology set $H^1(\G, \hpi_1(\bX,b))$ that classifies torsors, and a non-abelian Kummer map
$$\d^{na}: X(F) \ra H^1(\G, \hpi_1(\bX,b))$$
sending a point $x$ to the class of the torsor $\hpi_1(X;b,x)$. This is obviously a basic construction
whose importance, however, has begun to emerge only in the last twenty or so years. It relies very much on the flexible
use of varying base-points in Grothendieck's theory of the fundamental group, and it appears to have taken
some time after the inception of the arithmetic  $\pi_1$
theory \cite{SGA1} for the importance of such a variation to be properly appreciated
\cite{grothendieck, deligne, ihara}. In fact,
I believe the impetus for taking this
variation seriously came also for the most part from Hodge theory  \cite{hain1,
hain2}.
As far as Diophantine problems are concerned, in a letter to Faltings \cite{grothendieck}
written shortly after the proof of the Mordell conjecture,
Grothendieck proposed the remarkable conjecture
that  $\d^{na}$
 should be a bijection for compact hyperbolic curves. He expected such a statement to be directly relevant to the Mordell
problem and probably its variants like conjecture (2). This expectation appears still to be rather reasonable.
For one thing, it is evident that this conjecture is a hyperbolic analogue of the finiteness conjecture for
Tate-Shafarevich groups. And then, profound progress is represented by
 the work of Nakamura, Tamagawa, and
Mochizuki \cite{nakamura, tamagawa, mochizuki}, where a statement of this sort is proved when points
in the number field are replaced by
dominant maps from other varieties. Some marginal insight might also be gleaned from
\cite{kim1} and \cite{kim2} where a
unipotent analogue of the Kummer map is related to Diophantine finiteness theorems. There, the ambient space
inside which the analysis takes place is
 a classifying variety $H^1_f(\G_v, U^{et}_n)$ of torsors for the local unipotent \'etale
fundamental group (rather than the  Jacobian), while the finitely-generated
group $G$ is replaced by the image  of a map
$$H^1_f(\G_S, U^{et}_n)\ra H^1_f(\G_v, U^{et}_n)$$
coming from a space of global torsors. Thereby, one obtains a new manifestation of the analytic strategy
proving $X\cap Im[H^1_f(\G_S, U^{et}_n)]$
to be finite in some very special circumstances, and in general for a hyperbolic curve over $\Q$ if
one admits standard conjectures from the theory of mixed motives (for example, the Fontaine-Mazur conjecture
on geometric Galois representations). Fortunately, Chabauty's original method
 fits
naturally into this setting as the technical foundation of the
analytic part  now becomes non-abelian $p$-adic Hodge theory and iterated integrals.

It should be evident  at this point that the Galois theory of the title refers in general to the theory of the fundamental
group. Serge Lang was profoundly concerned with the fundamental group for a good part of his mathematical
life. A rather haphazard list of evidence might be comprised of:

-his foundational work on unramified class field theory for varieties over finite fields, where he proves the
surjectivity of the reciprocity map among many other things \cite{UCFT, SL};

-his study of the ubiquitous `Lang torsor' \cite{AG};

-his work with Serre on fundamental groups of
proper varieties in arbitrary characteristic \cite{LS};

-his extensive study with Kubert of the modular function field \cite{KuL};

-his work with Katz \cite{KL} on finiteness theorems for relative $\pi_1$'s
that made possible the subsequent proof by Bloch \cite{bloch}, and then Kato and Saito \cite{KS1,KS2} of the
finiteness of $CH_0$ for arithmetic schemes.

Besides these influential papers, the reader is referred to his beautiful AMS colloquium lectures \cite{UCG}
for a global perspective on the role of covering spaces in
arithmetic.

Even towards the end of his life when his published work
went in an increasingly  analytic direction, he
had a keen interest  both in fundamental groups and in
the analogy between hyperbolic manifolds and number fields wherein
fundamental
groups play a central role. In my last year of graduate school, he urged me strongly
to study the work of Kato and Saito (and apply it to Arakelov theory!)
even though it had been years since he had himself
been involved with such questions. From  the Spring of 2004, I recall a characteristically animated
exchange in the course of which
he explained to me a theorem of Geyer \cite{geyer}
stating that abelian subgroups of $\Gal
(\bQ/\Q)$ are pro-cyclic. It was  clear that he perceived  this fact to fit nicely into his vivid ideas about the heat kernel
\cite{JL}, but in
a manner that I failed (and still fail) to comprehend properly. (He was unfortunately
secretive with his deeper reflections on the arithmetic significance of his
later work, allowing only informal glimpses here and there. It is tempting but
probably
premature
to speculate about a Galois theory that  encompasses even
Archimedean analysis.)
The preoccupation
with hyperbolic geometry that was evident even from  the 70's (\cite{HDP, HDA, HS} and  \cite{NT}, chapters 8 and 9)
could rather generally be construed as reflecting a persistent intuition about the
relevance of fundamental groups to Diophantine problems. (An intuition that
was shared by Grothendieck \cite{grothendieck} and even Weil \cite{weil}.)

As for the direct application of non-abelian
fundamental groups to Diophantine geometry that we have outlined here, one
can convincingly place it into the general framework of Lang's inquiries.
He is discussing the theorem
of Siegel in the following paragraph from the notes to chapter 8 of \cite{FDG}:
\begin{quotation}
 The general version  used here was presented in \cite{IP} following Siegel's
(and Mahler's) method. The Jacobian replaces the theta function, as usual, and the mechanism of the covering
already used by Siegel appears here in its full formal clarity. It is striking to observe that in \cite{UCFT}, I used
the Jacobian in a formally analogous way to deal with the class field theory in function fields. In that case,
Artin's reciprocity law was reduced to a formal computation in the isogeny $u\mapsto u^{(q)}-u$ of the
Jacobian. In the present case, the heart of the proof is reduced to a formal computation of heights in the
isogeny $u\mapsto mu+a$.
\end{quotation}
We have emphasized above  the importance of the Kummer map
$$x \mapsto [\hpi_1(\bX;b,x)]\in H^1(\G, \hpi_1(\bX;b)).$$
When $X$ is defined over a finite field $\F_q$ and we replace $\hpi_1(\bX, b)$ by its abelian quotient
$H_1(\bX, \hZ)$, the map
takes values in $$H^1(\Gal (\bar{\F}_q/\F_q), H_1(\bX, \hZ))= H_1(\bX, \hZ)/[(Fr-1)H_1(\bX, \hZ)],$$
$Fr\in Gal (\bar{\F}_q/\F_q )$ being the
Frobenius element. But this last group is nothing but the kernel
$$\hpi_1^{ab}(X)^0$$ of the structure map
$$\hpi_1^{ab}(X) \ra \hpi_1(\Spec(\F_q)).$$
Thus the abelian quotient of the Kummer map becomes identified with the {\em reciprocity} map
\cite{KS1}
$$CH_0(X)^0 \ra \hpi_1^{ab}(X)^0$$
of unramified class field theory evaluated on the cycle
$(x)-(b)$. In other words, the reciprocity map is merely an `abelianized' Kummer map in this situation.
There is no choice but to interprete
the reciprocity law \cite{KS1, KS2} as an `abelianized Grothendieck conjecture' over finite fields.

 Of course it is hard to imagine
exactly what Lang himself found striking in the analogy when he wrote the lines quoted
above. What is not hard to
imagine is that he would have been very much at home with the ideas surrounding
Grothendieck's conjecture and the non-abelian Kummer map.

\end{document}